# On matrix algebras isomorphic to finite fields and planar Dembowski-Ostrom monomials

Christof Beierle [a,*], Patrick Felke [b]

[a] Ruhr University Bochum, Bochum, Germany
[b] University of Applied Sciences Emden/Leer, Emden, Germany

ABSTRACT

Let $p$ be a prime and $n$ a positive integer. As the first main result, we present a *deterministic* algorithm for deciding whether the matrix algebra $\mathbb{F}_p[A_1, \ldots, A_t]$ with $A_1, \ldots, A_t \in \mathrm{GL}(n, \mathbb{F}_p)$ is a finite field, performing at most $\mathcal{O}(tn^6 \log(p))$ elementary operations in $\mathbb{F}_p$. In the affirmative case, the algorithm returns a defining element $a$ so that $\mathbb{F}_p[A_1, \ldots, A_t] = \mathbb{F}_p[a]$.

We then study an invariant for the extended-affine equivalence of Dembowski-Ostrom (DO) polynomials. More precisely, for a DO polynomial $g \in \mathbb{F}_{p^n}[x]$, we associate to $g$ a set of $n \times n$ matrices with coefficients in $\mathbb{F}_p$, denoted $\mathrm{Quot}(\mathcal{D}_g)$, that stays invariant up to matrix similarity when applying extended-affine equivalence transformations to $g$. In the case where $g$ is a *planar* DO polynomial, $\mathrm{Quot}(\mathcal{D}_g)$ is the set of quotients $XY^{-1}$ with $Y \neq 0, X$ being elements from the spread set of the corresponding commutative presemifield, and $\mathrm{Quot}(\mathcal{D}_g)$ forms a field of order $p^n$ if and only if $g$ is equivalent to the planar monomial $x^2$, i.e., if and only if the commutative presemifield associated to $g$ is isotopic to a finite field.

As the second main result, we analyze the structure of $\mathrm{Quot}(\mathcal{D}_g)$ for all planar DO *monomials*, i.e., for commutative presemifields of odd order being isotopic to a finite field or a commutative twisted field. More precisely, for $g$ being equivalent to a planar DO monomial, we show that every non-zero element $X \in \mathrm{Quot}(\mathcal{D}_g)$ generates a field $\mathbb{F}_p[X] \subseteq \mathrm{Quot}(\mathcal{D}_g)$ and $\mathrm{Quot}(\mathcal{D}_g)$ *contains* the field $\mathbb{F}_{p^n}$.

* Corresponding author.
*E-mail addresses:* christof.beierle@rub.de (C. Beierle), patrick.felke@hs-emden-leer.de (P. Felke).

## 1. Introduction

A *semifield* is an algebraic generalization of a field in which the multiplication is neither required to be commutative nor associative. A *presemifield* does further not enforce the existence of a multiplicative identity. A (pre)semifield is called *commutative* if the multiplication is commutative. It is known that any finite (commutative) presemifield can be represented as a structure $\mathcal{R} = (\mathbb{F}_q, +, \star)$ where $\mathbb{F}_q$ is the finite field with $q$ elements, $(\mathbb{F}_q, +)$ the additive group of $\mathbb{F}_q$, and $\star$ a (commutative) binary operation on $\mathbb{F}_q$ such that $(\mathbb{F}_q \setminus \{0\}, \star)$ forms a quasigroup, $\mathbb{F}_q \star 0 = 0$, and $\star$ fulfills distributivity over $+$ (see [8,20]). Hence, the order of a finite presemifield is necessarily a prime power. In [8], Coulter and Henderson showed a one-to-one correspondence between commutative presemifields of odd order and planar Dembowski-Ostrom polynomials. We will recall their results below, after first introducing the required terminology.

Let $p$ be an odd prime and $n$ a positive integer. A polynomial $g \in \mathbb{F}_{p^n}[x]$ is called *planar* if, for all $\alpha \in \mathbb{F}_{p^n}^*$,

$$\Delta_{g,\alpha}(x) := g(x+\alpha) - g(x) - g(\alpha)$$

is a permutation polynomial in $\mathbb{F}_{p^n}[x]$ i.e., its evaluation map $\mathbb{F}_{p^n} \to \mathbb{F}_{p^n}, y \mapsto \Delta_{g,\alpha}(y)$ is 1-to-1. Planar polynomials were introduced by Dembowski and Ostrom in [10]. Since we only study properties of evaluation maps in $\mathbb{F}_{p^n}$, we assume that $g \in \mathbb{F}_{p^n}[x]/(x^{p^n} - x)$, i.e., $g$ has degree at most $p^n - 1$. A special type of polynomials in $\mathbb{F}_{p^n}[x]$ are *Dembowski-Ostrom* (DO) polynomials, which are those of the form

$$\sum_{0 \le i \le j \le n-1} u_{i,j} \cdot x^{p^i + p^j}, \quad u_{i,j} \in \mathbb{F}_{p^n}.$$

An equivalence relation between two polynomials that leaves the planarity property invariant is *graph equivalence* (also called *CCZ-equivalence*) [5]. The vast majority of the known classes of planar polynomials are equivalent to DO polynomials, and only one other family is known [9]. The graph equivalence of two planar polynomials coincides with *extended-affine equivalence*, and with *linear equivalence* in the case of planar DO polynomials [4].

Coulter and Henderson showed that any commutative presemifield $\mathcal{R} = (\mathbb{F}_q, +, \star)$ of odd order $q$ gives rise to a planar DO polynomial $g \in \mathbb{F}_q[x]$ via $g(a) = a \star a$ for all $a \in \mathbb{F}_q$, and conversely, any planar DO polynomial $g \in \mathbb{F}_q[x]$ defines a commutative presemifield $\mathcal{R}_g = (\mathbb{F}_q, +, \star)$ via $a \star b := \Delta_{g,a}(b)$. We further refer to Pott [27] for a survey on known classes of planar (DO) polynomials and their relation to commutative presemifields. When it comes to the problem of classifying (pre)semifields, an important concept is the notion of *isotopy*: Two finite presemifields $\mathcal{R} = (\mathbb{F}_q, +, \star)$ and $\mathcal{R}' = (\mathbb{F}_q, +, \circ)$ are called *isotopic* if there exist additive permutations $A, B, C$ over $\mathbb{F}_q$ such that the identity $A(a) \star B(b) = C(a \circ b)$ holds for all $a, b \in \mathbb{F}_q$. Any finite commutative presemifield $\mathcal{R} = (\mathbb{F}_q, +, \star)$ which is not already a semifield is isotopic to a commutative semifield by choosing an arbitrary $a \in \mathbb{F}_q^*$ and defining a new multiplication $\circ$ by $(x \star a) \circ (a \star y) = x \star y$ (see [8]). This is usually referred to as *Kaplansky's trick*. The problem of classifying commutative semifields (of even or odd order) is equivalent to the problem of classifying *symplectic semifields*, as every commutative semifield (up to isotopy) can be obtained within the Knuth orbit of a symplectic semifield and vice versa. We will omit the definition of a symplectic semifield and of the Knuth orbit in this article and instead refer to [18] and [21] for the details. Besides a classification of semifields of order $p^n$ with $n \leq 3$ ([20,25]), semifields are classified for order $2^n$ with $n \leq 6$, order $3^n$ with $n \leq 5$, order $5^4$, and order $7^4$ ([19,20,11,28,6,7,29]). Lavrauw and Sheekey [22] recently classified 4-dimensional (over $\mathbb{F}_8$ and $\mathbb{F}_9$, respectively) symplectic semifields of order $2^{12}$ and $3^8$. A complete classification of commutative (pre)semifields remains elusive and their number grows at least exponentially in their dimension [18,14].[1]

The authors of [8] showed a crucial property of those two latter classes of commutative semifields; namely, the commutative presemifield $\mathcal{R}_g$ is isotopic to a finite field or a commutative twisted field if and only if $g$ is equivalent to a planar DO *monomial*, i.e., a planar DO polynomial of the form $x^d$. In [9], Coulter and Matthews showed that any planar DO monomial in $\mathbb{F}_{p^n}[x]$ is equivalent to $x^{p^k+1} \in \mathbb{F}_{p^n}[x]$ with $n/\gcd(k, n)$ being odd. The case where $\mathcal{R}_g$ is isotopic to a finite field then corresponds to $g$ being equivalent to $x^2$ (see [8, Cor. 3.10]).

The state-of-the-art algorithm to decide the graph equivalence between two polynomials in $\mathbb{F}_{p^n}[x]$ relies on deciding code equivalence of two linear codes, see [13]. Recently, Ivkovic and Kaleyski developed a more efficient algorithm for deciding linear equivalence of $t$-to-1 functions [17]. Since the evaluation map of every planar DO polynomial is 2-to-1 (see [33]), this algorithm can be used to completely decide the equivalence between planar DO polynomials. Still, the complexity of this algorithm is *exponential* in $n$. As we will see, for deciding the equivalence of a DO polynomial to $x^2$, there are more efficient algorithms.

*Our results*

This work contains two main results. First, for a fixed prime $p$, we present a *deterministic* algorithm with time complexity polynomial in $n$ and $t$ for deciding whether an $\mathbb{F}_p$-algebra generated by $t$ matrices of dimension $n \times n$ is a finite field. In the affirmative case, the algorithm returns a defining element of the field extension over $\mathbb{F}_p$ and this element is computed deterministically.

[1] However, for $n$ being a prime all commutative semifields of order $p^n$ with $p$ being a large enough prime are isotopic to a finite field or a commutative twisted field [26].

**Theorem 1.** *Let $p$ be a prime, $n, t$ positive integers, and $\mathcal{S} = \{A_1, \dots, A_t\}$ be a subset of $\mathrm{GL}(n, \mathbb{F}_p)$. There exists a deterministic algorithm (viz., Algorithm 2) performing at most $\mathcal{O}(tn^6 \log(p))$ elementary operations in $\mathbb{F}_p$ and which decides whether $\mathbb{F}_p[\mathcal{S}]$ is a finite field. In the affirmative case, Algorithm 2 computes the extension degree $s$ over $\mathbb{F}_p$ and a defining element $a$ so that $\mathbb{F}_p[\mathcal{S}] = \mathbb{F}_p[a]$.*

A previously published algorithm ([3]) for this task[2] had to compute $\mathcal{O}(t \log p^n + tn^4 \log p \log \log p^n)$ elementary field operations and, most importantly, needed a factoring oracle for $p^n - 1$. Although it is well known that any finite field is a simple extension over its prime subfield, the standard textbook argument is non constructive and uses the cyclicity of the multiplicative group. However, if one is not interested in finding a primitive element of the multiplicative group, our algorithm can be used to *construct* a defining element of the field extension in polynomial time.

We then study an invariant for the extended-affine equivalence relation between two (not necessarily planar) DO polynomials. More precisely, we associate to a DO polynomial $g \in \mathbb{F}_{p^n}[x]$ a set of $n \times n$ matrices with coefficients in $\mathbb{F}_p$, denoted $\mathrm{Quot}(\mathcal{D}_g)$, that stays invariant up to the similarity transformation of matrices when applying extended-affine equivalence transformations to $g$ (Theorem 3 and Corollary 2). In the case where $g$ is a planar DO polynomial, $\mathrm{Quot}(\mathcal{D}_g)$ is the set of quotients $XY^{-1}$ with $X, Y \in \mathrm{GL}(n, \mathbb{F}_p) \cup \{0\}, Y \neq 0$ and $X, Y$ being elements from the spread set of the corresponding commutative presemifield. We recall that the *spread set* is the set of matrices corresponding to the mappings $x \mapsto a \star x$ of left-multiplications with elements $a$. This is a well-known invariant in the context of semifields and was already used in previous classification results, e.g., [28]. Since $\mathrm{Quot}(\mathcal{D}_g)$ is a finite field if and only if $g$ is equivalent to $x^2$, this allows the application of Algorithm 2 to decide the equivalence of a planar DO polynomial $g$ to $x^2$ (resp., to decide the isotopy of a commutative presemifield $\mathcal{R}_g$ to a finite field of odd order) using $\mathcal{O}(n^7 \log(p))$ elementary operations in $\mathbb{F}_p$ and $\mathcal{O}(n^2)$ evaluations of $g$ (resp., $\mathcal{O}(n^2)$ evaluations of $\star$ in $\mathcal{R}_g$).[3]

As our second main contribution, we analyze the structure of $\mathrm{Quot}(\mathcal{D}_g)$ for the remaining planar DO monomials, i.e., for those commutative presemifields of odd order being isotopic to a commutative twisted field [1]. More precisely, for $g \in \mathbb{F}_{p^n}[x]$ equivalent to $x^{p^k+1}$ with $n/\gcd(k, n)$ being odd, we show that any non-zero element in $\mathrm{Quot}(\mathcal{D}_g)$ generates a field *contained* in $\mathrm{Quot}(\mathcal{D}_g)$ and, in particular, $\mathrm{Quot}(\mathcal{D}_g)$ always contains a finite field of order $p^n$. Denoting the matrix of the mapping $x \mapsto a \star x$ in the corresponding commutative presemifield by $M_{g,a}$, our main theorem is stated as follows.

[2] Compared to our approach, the algorithm in [3] returns a generator of the multiplicative group $\langle A_1, \dots, A_t \rangle$ in case $\mathbb{F}_p[A_1, \dots, A_t]$ is a field.

[3] However, as this particular problem could also be solved more efficiently by means of the nuclei of $\mathcal{R}_g$ (Remark 4), this application of Algorithm 2 is not one of our main contributions.

**Theorem 2.** *Let $p$ be an odd prime and $n$ a positive integer. Let $g(x) \in \mathbb{F}_{p^n}[x]$ be a planar DO monomial. For any $\alpha, \beta \in \mathbb{F}_{p^n}^*$, the element $X := M_{g,\beta} M_{g,\alpha}^{-1} \in \mathrm{Quot}(\mathcal{D}_g)$ generates a field isomorphic to $\mathbb{F}_p(\alpha^{-1}\beta)$ viz. $\mathbb{F}_p[X]$, and $\mathbb{F}_p[X] \subseteq \mathrm{Quot}(\mathcal{D}_g)$.*

While it was known before that any element in $\mathrm{Quot}(\mathcal{D}_g)$ generates a finite field, the crucial part of our result is the *containment* of this field within $\mathrm{Quot}(\mathcal{D}_g)$ itself.

## 2. Preliminaries

By $\mathrm{Mat}_{\mathbb{F}_p}(n, n)$, we denote the ring of all $n \times n$ matrices with coefficients in the prime field $\mathbb{F}_p$ and by $\mathrm{GL}(n, \mathbb{F}_p)$ the group of all invertible matrices in $\mathrm{Mat}_{\mathbb{F}_p}(n, n)$. For $A \in \mathrm{Mat}_{\mathbb{F}_p}(n, n)$, we denote its minimal polynomial by $\mu_A$. Given a non-empty set $\mathcal{S} = \{A_1, \ldots, A_t\} \subseteq \mathrm{Mat}_{\mathbb{F}_p}(n, n)$, we denote by $\mathbb{F}_p[\mathcal{S}]$ (or $\mathbb{F}_p[A_1, \ldots, A_t]$) the $\mathbb{F}_p$-algebra generated by $\mathcal{S}$, i.e., the intersection of all subalgebras of $\mathrm{Mat}_{\mathbb{F}_p}(n, n)$ containing $\mathcal{S}$. Any finite field $\mathbb{F}_{p^n}$ (resp., a proper subfield $\mathbb{F}_{p^m}$) is isomorphic to $\mathbb{F}_p[T_\beta]$, where $T_\beta$ denotes a matrix corresponding to the linear mapping $x \mapsto \beta x$ over $\mathbb{F}_{p^n}$, for $\beta \in \mathbb{F}_{p^n}^*$ defining a polynomial basis of $\mathbb{F}_{p^n}$ (resp., of $\mathbb{F}_{p^m}$). For more details on *matrix representations* of finite fields, we refer to, e.g., [15] or [23]. The following statement is well known, see e.g., [15, Ch. 7.2].

**Lemma 1.** *Let $\mathcal{S}, \mathcal{T} \subseteq \mathrm{GL}(n, \mathbb{F}_p)$ be non-empty and let $\mathbb{F}_p[\mathcal{S}]$ be a matrix algebra isomorphic to $\mathbb{F}_{p^n}$. The matrix algebra $\mathbb{F}_p[\mathcal{T}]$ is isomorphic to $\mathbb{F}_{p^n}$ if and only if there exists an $A \in \mathrm{GL}(n, \mathbb{F}_p)$ such that $\mathbb{F}_p[\mathcal{T}] = A^{-1} \cdot \mathbb{F}_p[\mathcal{S}] \cdot A$.*

A polynomial $L(x) \in \mathbb{F}_{p^n}[x]$ of the form

$$L(x) = \sum_{i=0}^{m} u_i x^{p^i}, \quad u_i \in \mathbb{F}_{p^n}$$

is called a *linearized polynomial*. Given linearized polynomials $L, L' \in \mathbb{F}_{p^n}[x]$, their *symbolic product*, denoted $\otimes$, is defined as their composition, i.e.,

$$L(x) \otimes L'(x) := L(L'(x)).$$

Every linearized polynomial $L \in \mathbb{F}_{p^n}[x]$ induces a linear mapping over $\mathbb{F}_{p^n}$ via its evaluation map $y \mapsto L(y)$. Vice versa, for every linear mapping over $\mathbb{F}_{p^n}$ there exists a unique corresponding linearized polynomial in $\mathbb{F}_{p^n}[x]/(x^{p^n} - x)$ of degree at most $p^n - 1$. Hence, linearized polynomials in $\mathbb{F}_{p^n}[x]/(x^{p^n} - x)$ are closed under taking symbolic products and reducing modulo $x^{p^n} - x$. By fixing a basis $B$ of $\mathbb{F}_{p^n}$ as an $\mathbb{F}_p$-vector space, each linear mapping $\phi$ over $\mathbb{F}_{p^n}$ can be associated to a matrix $M \in \mathrm{Mat}_{\mathbb{F}_p}(n, n)$ such that $\phi(x)$ corresponds to the matrix multiplication $M \cdot x$ in $\mathbb{F}_p^n$. Therefore, (after fixing the basis $B$) we can associate to $L$ a matrix $M \in \mathrm{Mat}_{\mathbb{F}_p}(n, n)$ such that the linear mapping induced by $L$ corresponds to $M$. The link between a linearized polynomial in $\mathbb{F}_{p^n}[x]$ and an element in $\mathrm{Mat}_{\mathbb{F}_p}(n,n)$ is made explicit in, e.g., [35]. We have the following fact.

**Lemma 2.** *Let $M, M'$ be the matrices associated to the linearized polynomials $L$ and $L'$ in $\mathbb{F}_{p^n}[x]$, respectively. Then, the matrix associated to $L \otimes L' \mod (x^{p^n} - x)$ is equal to $M \cdot M'$, and the matrix associated to $L + L'$ is equal to $M + M'$.*

In the following, we always assume a fixed choice for the basis $B$, so that we can switch between linearized polynomials and matrices. For a DO polynomial $g(x) \in \mathbb{F}_{p^n}[x]$, we denote by $\Delta_{g,\alpha}$ the *linearized derivative of $g$ in direction $\alpha$*, i.e., the polynomial $g(x+\alpha) - g(x) - g(\alpha) \in \mathbb{F}_{p^n}[x]$. Since $g$ is DO, $\Delta_{g,\alpha}$ is a linearized polynomial. Let us denote by $M_{g,\alpha}$ the matrix in $\mathrm{Mat}_{\mathbb{F}_p}(n,n)$ associated to $\Delta_{g,\alpha}$. We define the set $\mathcal{D}_g$ of *linearized derivative matrices* as

$$\mathcal{D}_g := \{M_{g,\alpha} \mid \alpha \in \mathbb{F}_{p^n}\}.$$

If $g$ is planar, the set $\mathcal{D}_g$ is equal to the set of matrices corresponding to the mappings $x \mapsto \alpha \star x$ of left-multiplications with elements $\alpha$ in the corresponding commutative presemifield $\mathcal{R}_g$. In that case, $\mathcal{D}_g$ is also called the *spread set of $\mathcal{R}_g$* (see, e.g., [30, Sec. 2.1]).

Two polynomials $g, g' \in \mathbb{F}_{p^n}[x]/(x^{p^n} - x)$ are called *extended-affine equivalent* (which we shortly call *equivalent*), if there exist linearized permutation polynomials $L, L' \in \mathbb{F}_{p^n}[x]$, a linearized polynomial $L'' \in \mathbb{F}_{p^n}[x]$, and constants $u, v \in \mathbb{F}_{p^n}$ such that

$$g'(x) = L'(g(L(x)+u)) + v + L''(x) \mod (x^{p^n} - x).$$

For the special case of $L'' = 0$ and $u = v = 0$, the polynomials $g(x)$ and $g'(x)$ are called *linear equivalent.* For planar DO polynomials, we know that extended-affine equivalence (even the more general CCZ-equivalence) coincides with linear equivalence [4].

## 3. Deciding FiniteField in polynomial time

Let $\mathcal{S} = \{A_1, \ldots, A_t\}$ be a subset of $\mathrm{GL}(n, \mathbb{F}_p)$. Then, $\mathbb{F}_p[\mathcal{S}]$ is a field if and only if the multiplicative group generated by $\mathcal{S}$ is cyclic and generated by an element with irreducible minimal polynomial (see [3, Thm. 1]). This can be decided (together with constructing a generator of $\langle \mathcal{S} \rangle$ in case $\mathbb{F}_p[\mathcal{S}]$ is a field) by a deterministic algorithm which performs $\mathcal{O}(t \log p^n + tn^4 \log p \log\log p^n)$ elementary field operations in $\mathbb{F}_p$, supposing the prime factorization of $p^n - 1$ is known, see [3, Thm. 2]. The knowledge of the prime factorization of $p^n - 1$ is required as the algorithm needs to compute several multiplicative orders of elements in $\mathrm{GL}(n, \mathbb{F}_p)$. However, if we just want to find a defining element of the field extension over its prime subfield *and do not need to find a generator of $\langle \mathcal{S} \rangle$*, we can do so by a deterministic polynomial-time (with a constant $\log(p)$) algorithm which does not need an oracle for the prime factorization of $p^n - 1$. More precisely, we will prove the following theorem.

**Theorem 1.** *Let $p$ be a prime, $n, t$ positive integers, and $\mathcal{S} = \{A_1, \ldots, A_t\}$ be a subset of $\mathrm{GL}(n, \mathbb{F}_p)$. There exists a deterministic algorithm (viz., Algorithm 2) performing at most $\mathcal{O}(tn^6 \log(p))$ elementary operations in $\mathbb{F}_p$ and which decides whether $\mathbb{F}_p[\mathcal{S}]$ is a finite field. In the affirmative case, Algorithm 2 computes the extension degree $s$ over $\mathbb{F}_p$ and a defining element $a$ so that $\mathbb{F}_p[\mathcal{S}] = \mathbb{F}_p[a]$.*

The proof requires Lemma 3, which will be given after we have fixed some notation. For a ring element $x$ and positive integers $k, \ell$ with $\ell$ being a multiple of $k$, we define $\mathrm{Tr}_{\ell,k}(x) := x + x^q + x^{q^2} + \cdots + x^{q^{\frac{\ell}{k}-1}}$, where $q := p^k$. Recall that when taking $x$ from a finite field $\mathbb{F}_{p^\ell}$, the mapping $\mathrm{Tr}_{\ell,k}$ is the *relative trace* mapping, which is an $\mathbb{F}_{p^k}$-linear mapping from $\mathbb{F}_{p^\ell}$ onto $\mathbb{F}_{p^k}$. Given a finite field extension $L/K$, we denote by $[L : K]$ its extension degree, i.e., the dimension of $L$ as a $K$-vector space. By $\overline{K}$, we denote the algebraic closure of $K$, and for $a_1, a_2, \ldots, a_t \in \overline{K}$, we denote by $K(a_1, \ldots, a_t)$ the field of adjoining $a_1, \ldots, a_t$ to $K$, i.e., the smallest field containing $K$ and $a_1, \ldots, a_t$. It is well known that any finite extension of a finite field $K$ is simple, i.e., it is generated by adjoining only one element to $K$. The usual argument for that fact given in the literature is using the cyclicity of the multiplicative group of a finite field. In our case, we are interested in efficiently *computing* such a defining element. The following lemma provides us a way to do so. Note that, while (1) is a trivial fact, (2) was already pointed out in previous literature ([24, Rem. 2.3]).

**Lemma 3.** *Given $m \le n$ and $a, b \in \overline{\mathbb{F}}_p$ such that $[\mathbb{F}_p(a) : \mathbb{F}_p] = m$ and $[\mathbb{F}_p(b) : \mathbb{F}_p] = n$, i.e. $\mathbb{F}_p(a) = \mathbb{F}_{p^m}$ and $\mathbb{F}_p(b) = \mathbb{F}_{p^n}$. Then, the following assertions hold.*

(1) *If $m$ divides $n$, we have $\mathbb{F}_p(a, b) = \mathbb{F}_p(b) = \mathbb{F}_{p^{\mathrm{lcm}(m,n)}} = \mathbb{F}_{p^n}$.*
(2) *If $\gcd(m, n) = 1$, we have $\mathbb{F}_p(a, b) = \mathbb{F}_p(ab) = \mathbb{F}_{p^{\mathrm{lcm}(m,n)}} = \mathbb{F}_{p^{mn}}$.*
(3) *Let $m = \kappa_1^{e_{m,1}} \cdot \kappa_2^{e_{m,2}} \cdot \cdots \cdot \kappa_r^{e_{m,r}}$ and $n = \kappa_1^{e_{n,1}} \cdot \kappa_2^{e_{n,2}} \cdot \cdots \cdot \kappa_r^{e_{n,r}}$ for pairwise distinct primes $\kappa_i$ and non-negative integers $e_{m,i}, e_{n,i}$. Then, for each $i = 1, \ldots, r$, there exists $\ell_{m,i} \in \{1, \ldots, m-1\}$ and $\ell_{n,i} \in \{1, \ldots, n-1\}$ such that $[\mathbb{F}_p(a^{(i)}) : \mathbb{F}_p] = \kappa_i^{e_{m,i}}$, $[\mathbb{F}_p(b^{(i)}) : \mathbb{F}_p] = \kappa_i^{e_{n,i}}$, where $a^{(i)} := \mathrm{Tr}_{m,\kappa_i^{e_{m,i}}}(a^{\ell_{m,i}})$ and $b^{(i)} := \mathrm{Tr}_{n,\kappa_i^{e_{n,i}}}(b^{\ell_{n,i}})$. For*

$$c^{(i)} := \begin{cases} a^{(i)} & \text{if } e_{m,i} \ge e_{n,i} \\ b^{(i)} & \text{otherwise,} \end{cases}$$

*we then have $\mathbb{F}_p(a, b) = \mathbb{F}_p(\prod_{i=1}^r c^{(i)}) = \mathbb{F}_{p^{\mathrm{lcm}(m,n)}}$.*

**Proof.** Assertion 1 is trivial since, if $m$ divides $n$, we have $\mathbb{F}_p(a) \subseteq \mathbb{F}_p(b)$ and obviously $\mathbb{F}_{p^n} = \mathbb{F}_p(b) = \mathbb{F}_p(a, b) = \mathbb{F}_{p^{\mathrm{lcm}(m,n)}}$.

Assertion 2 along with a proof is given in [24, Rem. 2.3]. We outline the proof (by contradiction) for completeness. Let us assume that $\mathbb{F}_p(ab)$ is a proper subfield of $\mathbb{F}_p(a,b)$. Then, $l := [\mathbb{F}_p(a,b) : \mathbb{F}_p(ab)] > 1$. Recall that $\mathbb{F}_p(a,b)$ is a Galois extension of $\mathbb{F}_p(ab)$ with cyclic Galois group $G$ of order $l$. Since $l > 1$, there exists a prime divisor $d$ of $l$ and thus a subgroup $H$ of $G$ of order $d$. Then, $H$ corresponds to a unique fixed field $F$ with $[\mathbb{F}_p(a,b) : F] = d$. By assumption, $d$ cannot divide both $m$ and $n$. Without loss of generality, let us assume that $d \nmid m$. Then, $a \in F$ as otherwise, $[F(a) : F] = d$ as $d$ is prime and $F(a) = \mathbb{F}_p(a,b)$, which would imply that $d$ divides $[\mathbb{F}_p(a) : \mathbb{F}_p] = m$. With $a, ab \in F$, also $b \in F$ and consequently $F = \mathbb{F}_p(a,b)$, a contradiction.

To prove Assertion 3, let us assume we have $m = \kappa_1^{e_{m,1}} \cdot \kappa_2^{e_{m,2}} \cdot \dots \cdot \kappa_r^{e_{m,r}}$ for pairwise distinct primes $\kappa_i$ and non-negative integers $e_{m,i}$. Let us fix an index $i \in \{1,\dots,r\}$. As $\mathrm{Tr}_{m,\kappa_i^{e_{m,i}}}$ is $\mathbb{F}_p$-linear, maps onto $\mathbb{F}_{p^{(\kappa_i^{e_{m,i}})}}$ and $\{1, a, a^2, \dots, a^{m-1}\}$ is an $\mathbb{F}_p$-basis of $\mathbb{F}_{p^m}$, there exists an exponent $\ell_{m,i} \in \{1,\dots,m-1\}$ such that $\mathrm{Tr}_{m,\kappa_i^{e_{m,i}}}(a^{\ell_{m,i}})$ lies not in a proper subfield of $\mathbb{F}_{p^{(\kappa_i^{e_{m,i}})}}$. Hence, $[\mathbb{F}_p(\mathrm{Tr}_{m,\kappa_i^{e_{m,i}}}(a^{\ell_{m,i}})) : \mathbb{F}_p] = \kappa_i^{e_{m,i}}$. The same argument holds for the statement on $n$. By construction, for any two distinct $i,j \in \{1,\dots,r\}$, the extension degrees $[\mathbb{F}_p(c^{(i)}) : \mathbb{F}_p]$ and $[\mathbb{F}_p(c^{(j)}) : \mathbb{F}_p]$ are coprime, so we have $\mathbb{F}_p(c^{(1)}, c^{(2)}, \dots, c^{(r)}) = \mathbb{F}_p(\prod_{i=1}^r c^{(i)}) = \mathbb{F}_{p^{\mathrm{lcm}(m,n)}} = \mathbb{F}_p(a,b)$ by Assertion 2. □

**Corollary 1.** *Let $p$ be a prime and $n$ a positive integer. Algorithm 1 runs in at most $\mathcal{O}(n^6 \log(p))$ elementary operations in $\mathbb{F}_p$. Given as input $a, b \in \mathbb{F}_{p^n}$, it computes $c \in \mathbb{F}_{p^n}$ with $\mathbb{F}_p(a,b) = \mathbb{F}_p(c)$.*

**Proof.** For any $\gamma \in \mathbb{F}_{p^n}$, the extension degree of the simple extension $\mathbb{F}_p(\gamma)/\mathbb{F}_p$ is given by the degree of the minimal polynomial $\mu_\gamma \in \mathbb{F}_p[x]$ of $\gamma$. Note that $\mu_\gamma$ is equal to the minimal polynomial of the linear mapping $x \mapsto \gamma x$ over $\mathbb{F}_{p^n}$. The correctness then follows immediately from Lemma 3.

The algorithm performs at most $\mathcal{O}(r \cdot \max(m,k))$ computations of a minimal polynomial of an element in $\mathrm{GL}(n, \mathbb{F}_p)$ and at most $\mathcal{O}(r \cdot \max(m,k))$ computations of a relative trace of a power of an element in $\mathrm{GL}(n, \mathbb{F}_p)$, where the exponent is bounded by $\max(m-1, k-1)$. The minimal polynomial of an element in $\mathrm{GL}(n, \mathbb{F}_p)$ can be computed in $\mathcal{O}(n^3)$ elementary field operations [32]. Raising a matrix $a \in \mathrm{GL}(n, \mathbb{F}_p)$ to a power $j$ can be done with repeated squaring in $\mathcal{O}(n^3 \log(j))$ operations. Any of the relative trace maps of $A := a^j$ can be evaluated by computing the $n$ powers $A, A^p, \dots, A^{p^{n-1}}$ and performing at most $n-1$ matrix additions. Using the repeated squaring algorithm, the complexity to do so is $\mathcal{O}(n \cdot n^3 \log(p))$. Since $r \le \max(m,k) \le n$, we obtain the complexity as claimed. □

Using Algorithm 1 as a subroutine, we can now decide in polynomial time whether a matrix algebra with $t$ given generators is a finite field.

**Proof of Theorem 1.** We show that Algorithm 2 fulfills the claims in the statement of the theorem. We first show its correctness: Suppose that $\mathbb{F}_p[\mathcal{S}] \subseteq \mathbb{F}_{p^n}$ is a finite field.

**Algorithm 1** ComputeGenerator.

**Input:** Elements $a, b \in \mathrm{GL}(n, \mathbb{F}_p)$.
**Output:** If $\mathbb{F}_p[a, b]$ is a field, an element $c \in \mathrm{GL}(n, \mathbb{F}_p)$ such that $\mathbb{F}_p[a, b] = \mathbb{F}_p[c]$.

$m \leftarrow \deg(\mu_a), \quad k \leftarrow \deg(\mu_b), \quad d \leftarrow \gcd(m, k)$
**if** $d = k$ **then** ▷ See 1. of Lemma 3
  **return** $a$
**end if**
**if** $d = m$ **then** ▷ See 1. of Lemma 3
  **return** $b$
**end if**
**if** $d = 1$ **then** ▷ See 2. of Lemma 3
  **return** $ab$
**end if**
compute a list $[(q_i, e_{m,i}, e_{k,i})]_{i=1}^{r}$ s.t. $m = q_1^{e_{m,1}} \cdots q_r^{e_{m,r}}$ and $k = q_1^{e_{k,1}} \cdots q_r^{e_{k,r}}$ for pairwise distinct primes $q_i$ and $e_{m,i}, e_{k,i} \in \mathbb{N} \cup \{0\}$, $r \leq \max(m, k)$
**for** $i = 1, \dots, r$ **do**
  **if** $e_{m,i} \geq e_{k,i}$ **then**
    **for** $j = 1, \dots, m-1$ **do**
      $c^{(i)} \leftarrow \mathrm{Tr}_{m, q_i^{e_{m,i}}}(a^j)$
      **if** $\deg(\mu_{c^{(i)}}) = q_i^{e_{m,i}}$ **then**
        **break**
      **end if**
    **end for**
  **end if**
  **if** $e_{m,i} < e_{k,i}$ **then**
    **for** $j = 1, \dots, k-1$ **do**
      $c^{(i)} \leftarrow \mathrm{Tr}_{k, q_i^{e_{k,i}}}(b^j)$
      **if** $\deg(\mu_{c^{(i)}}) = q_i^{e_{k,i}}$ **then**
        **break**
      **end if**
    **end for**
  **end if**
**end for**
**return** $\prod_{i=1}^{r} c^{(i)}$ ▷ See 3. of Lemma 3

Then, by Lemma 3, the element $a$ computed after Line 4 of Algorithm 2 is a generator of $\mathbb{F}_p[\mathcal{S}]$ as a field, i.e., $\mathbb{F}_p[\mathcal{S}] = \mathbb{F}_p[a]$. Necessarily, the minimal polynomial $\mu_a$ of $a$ must be irreducible (see, e.g., [3, Lem. 1]). Further, all elements $A_1, \dots, A_t$ must be in the linear span of $(1, a, a^2, \dots, a^{n-1})$ as $a$ defines a polynomial basis of $\mathbb{F}_p[a]$. The extension degree of $\mathbb{F}_p[a]$ over $\mathbb{F}_p$ is given by the degree of $\mu_a$, which is the output of the algorithm. Conversely, suppose that $\mathbb{F}_p[\mathcal{S}]$ is not a field. If Algorithm 2 does not output false in Line 6, the matrix algebra $\mathbb{F}_p[a]$ computed in Algorithm 2 is a field since $\mu_a$ is irreducible. If the algorithm further does not output false in Line 10, we have $A_1, \dots, A_t \in \mathbb{F}_p[a]$, hence $\mathbb{F}_p[\mathcal{S}] \subseteq \mathbb{F}_p[a]$. Since any subring of a finite field is a field, we obtain a contradiction to the assumption that $\mathbb{F}_p[\mathcal{S}]$ is not a field.

To show the bound on the complexity, we observe that Algorithm 1 is called $t-1$ times as a subroutine, hence we obtain $\mathcal{O}(tn^6 \log(p))$ elementary operations in $\mathbb{F}_p$ as the complexity until Line 4. The complexity of the remaining steps can be neglected. Indeed, computing the minimal polynomial of $a$ and deciding its irreducibility can be performed in $\mathcal{O}(n^3)$ and $\mathcal{O}(n^3 \log(p))$ elementary field operations, respectively (see [32], resp., [31, Thm. 20.1]). In Line 9 (which is performed $t$ times), we only need to solve a linear system with $n^2$ equations and $n$ unknowns over $\mathbb{F}_p$. $\square$

**Algorithm 2** FiniteField.

**Input:** Matrices $A_1, \dots, A_t \in \mathrm{GL}(n, \mathbb{F}_p)$.
**Output:** The extension degree $[\mathbb{F}_p[A_1, \dots, A_t] : \mathbb{F}_p]$ and a defining element if $\mathbb{F}_p[A_1, \dots, A_t]$ is a field, **false** otherwise.

1: $a \leftarrow A_1$
2: **for** $i = 2, \dots, t$ **do**
3: $\quad a \leftarrow$ComputeGenerator$(a, A_i)$
4: **end for**
5: **if** $\mu_a$ is not irreducible **then** ▷ $\mathbb{F}_p[a]$ is not a field
6: $\quad$ **return false**
7: **end if**
8: **for** $i = 1, \dots, t$ **do** ▷ Check whether $A_1, \dots, A_t$ are elements of the field $\mathbb{F}_p[a]$
9: $\quad$ **if** $A_i \notin \mathrm{Span}(1, a, a^2, \dots, a^{n-1})$ **then**
10: $\quad\quad$ **return false**
11: $\quad$ **end if**
12: **end for**
13: **return** $(a, \deg(\mu_a))$

The predecessor of Algorithm 2 given in [3] (which needs a factoring oracle of $p^n - 1$, but also yields a generator of the multiplicative group $\langle A_1, \dots, A_t \rangle$ in case $\mathbb{F}_p[A_1, \dots, A_t]$ is a field) was applied to identify structure induced by extension fields within the linear layer of a block cipher. More precisely, one of the problems studied in [3] was to decide whether a given $ms \times ms$ matrix with entries in a finite field $\mathbb{F}_p$ can be represented as an $m \times m$ matrix over the extension field $\mathbb{F}_{p^s}$. Of course, we can also use Algorithm 2 to solve this problem. Another application of Algorithm 2 will be given in Section 5.

**Remark 1.** Of course, there is a straightforward probabilistic method for deciding whether $\mathbb{F}_p[A_1, \dots, A_t]$ is a finite field (and returning a defining element in the affirmative case), by randomly choosing elements $a$ from $\mathbb{F}_p[A_1, \dots, A_t]$ and conducting the check in Lines 8–12 in Algorithm 2, so that the call of ComputeGenerator can be avoided. Indeed, in a finite field with $p^s$ elements, all elements not lying in a proper subfield can be chosen as defining elements, so most elements would succeed the check in Lines 8–12. More precisely, this number is equal to $s \cdot N_p(s)$, where $N_p(s)$ is the number of monic irreducible polynomials of degree $s$ over $\mathbb{F}_p$, i.e.,

$$s \cdot N_p(s) = \sum_{d|s} \mu(d) \cdot p^{\frac{s}{d}},$$

where $\mu$ is the Möbius function (see, e.g., [23]). However, this probabilistic approach does not decide the problem as there is the possibility of always picking elements lying in a proper subfield.

## 4. The quotients of linearized derivative matrices

Let $g \in \mathbb{F}_{p^n}[x]$ be a DO polynomial. We study the *set of quotients in* $\mathcal{D}_g$, defined as

$$\mathrm{Quot}(\mathcal{D}_g) := \bigcup_{Y \in \mathcal{D}_g \cap \mathrm{GL}(n, \mathbb{F}_p)} \mathcal{D}_g Y^{-1} = \{XY^{-1} \mid X, Y \in \mathcal{D}_g \text{ and } Y \text{ is invertible}\}.$$

The following observation is immediate from the fact that $g(x+y) - g(x) - g(y)$ is symmetric in $x$ and $y$ and bilinear.

**Lemma 4.** *Let $g \in \mathbb{F}_{p^n}[x]$ be a DO polynomial and $\{\alpha_1, \alpha_2, \ldots, \alpha_n\}$ be an $\mathbb{F}_p$-basis of $\mathbb{F}_{p^n}$. For each $Y \in \mathrm{GL}(n, \mathbb{F}_p)$, the set $\mathcal{D}_g Y^{-1}$ is an $\mathbb{F}_p$-vector space spanned by*

$$\{M_{g,\alpha_1} Y^{-1}, M_{g,\alpha_2} Y^{-1}, \ldots, M_{g,\alpha_n} Y^{-1}\}.$$

*Moreover, if $g$ is planar, the space $\mathcal{D}_g Y^{-1}$ is $n$-dimensional over $\mathbb{F}_p$.*

We remark that the set of quotients in $\mathcal{D}_g$ could be empty. In the general case, the following can be deduced on the maximal size of $\mathrm{Quot}(\mathcal{D}_g)$.

**Lemma 5.** *For any DO polynomial $g(x) \in \mathbb{F}_{p^n}[x]$, an upper bound on $|\mathrm{Quot}(\mathcal{D}_g)|$ is given by*

$$\frac{(p^n - p) \cdot (p^n - 1)}{p - 1} + p.$$

**Proof.** Let $g \in \mathbb{F}_{p^n}[x]$ be a DO polynomial. For any $\alpha \in \mathbb{F}_{p^n}^*, c \in \mathbb{F}_p^*$, we have $M_{g,c\alpha} = c \cdot M_{g,\alpha}$. For the inverse (if it exists), this implies $M_{g,c\alpha}^{-1} = c^{-1} \cdot M_{g,\alpha}^{-1}$. Thus, for any $\alpha, \beta \in \mathbb{F}_{p^n}^*$ and $c, d \in \mathbb{F}_p^*$ such that $M_{g,\beta}$ is invertible, we have $M_{g,c\alpha} M_{g,d\beta}^{-1} = cd^{-1} \cdot M_{g,\alpha} M_{g,\beta}^{-1}$. Let $B := |\mathcal{D}_g \cap \mathrm{GL}(n, \mathbb{F}_p)|$. As a multiset, $\mathrm{Quot}(\mathcal{D}_g) \setminus \{cI_n \mid c \in \mathbb{F}_p\}$ contains $(p^n - p)B$ elements and the multiplicity of each element is at least $p - 1$. Hence, as a set, we have

$$|\mathrm{Quot}(\mathcal{D}_g) \setminus \{cI_n \mid c \in \mathbb{F}_p\}| \leq \frac{(p^n - p)B}{p - 1}$$

and the result follows because $B \leq p^n - 1$. □

**Remark 2.** If, for a given DO polynomial $g(x) \in \mathbb{F}_{p^n}[x]$, the set $\mathrm{Quot}(\mathcal{D}_g)$ is not empty, it must contain the field $\mathbb{F}_p$ as a subset, viz., $\{M_{g,c\alpha} M_{g,\alpha}^{-1} \mid c \in \mathbb{F}_p\}$ for $M_{g,\alpha} \in \mathcal{D}_g$ being invertible. This is because the existence of an invertible element $M_{g,\alpha} \in \mathcal{D}_g$ implies $M_{g,c\alpha} M_{g,\alpha}^{-1} = c M_{g,\alpha} M_{g,\alpha}^{-1} = cI_n \in \mathrm{Quot}(\mathcal{D}_g)$ for all $c \in \mathbb{F}_p$. As we will show later, for the special case where $g$ is equivalent to a planar DO monomial, $\mathrm{Quot}(\mathcal{D}_g)$ contains the extension field $\mathbb{F}_{p^n}$ as a subset.

There exist planar DO polynomials attaining the upper bound given in Lemma 5 with equality, e.g., $x^{p+1} \in \mathbb{F}_{3^5}[x]$. However, the set $\mathrm{Quot}(\mathcal{D}_g)$ can be much smaller. Notice that for a planar DO polynomial $g(x) \in \mathbb{F}_{p^n}[x]$, the set $\mathrm{Quot}(\mathcal{D}_g)$ must contain at least $p^n$ elements. This lower bound on $|\mathrm{Quot}(\mathcal{D}_g)|$ is tight, as it is attained by equality for the planar function $g(x) = x^2$ (see Lemma 7).

The set of quotients in $\mathcal{D}_g$ can be used to derive an invariant for extended-affine equivalence of DO polynomials. The following lemma describes how the linearized derivative changes under applying extended-affine equivalence transformations.

**Lemma 6.** *Let $g(x)$ and $g'(x)$ be two DO polynomials in $\mathbb{F}_{p^n}[x]$ that are extended-affine equivalent via $g'(x) = L'(g(L(x)+u)) + v + L''(x) \mod (x^{p^n} - x)$, where $L, L'$ are linearized permutation polynomials and $L''$ is a linearized polynomial in $\mathbb{F}_{p^n}[x]$ and $u, v \in \mathbb{F}_{p^n}$. For all $\alpha \in \mathbb{F}_{p^n}$, we then have*

$$\Delta_{g',\alpha}(x) = L'(x) \otimes \Delta_{g,L(\alpha)}(x) \otimes L(x) \mod (x^{p^n} - x). \tag{1}$$

**Proof.** Let $\alpha, y \in \mathbb{F}_{p^n}$. Given $g'(y) = L'(g(L(y)+u)) + v + L''(y)$, we can verify that

$$\Delta_{g',\alpha}(y) = L'\left(g(L(y+\alpha)+u) - g(L(y)+u) - g(L(\alpha)+u)\right) - v.$$

Substituting $\alpha$ by $L^{-1}(\alpha)$ yields

$$\Delta_{g',L^{-1}(\alpha)}(y) = L'\left(g(L(y)+\alpha+u) - g(L(y)+u) - g(\alpha+u)\right) - v,$$

thus

$$v + \Delta_{g',L^{-1}(\alpha)}(L^{-1}(y-u)) = L'\left(g(y+\alpha) - g(y) - g(\alpha+u)\right)$$
$$= L'\left(\Delta_{g,\alpha}(y) - \Delta_{g,\alpha}(u) - g(u)\right) = L'(\Delta_{g,\alpha}(y-u)) - L'(g(u)),$$

where the last equality follows from the linearity of $y \mapsto \Delta_{g,\alpha}(y)$. Substituting $y-u$ by $y$ yields the desired equality with $v + L'(g(u)) = 0$ because both $g$ and $g'$ are DO. □

When formulating Equation (1) in the language of matrices, we obtain

$$M_{g',\alpha} = M_{L'} \cdot M_{g,L(\alpha)} \cdot M_L, \tag{2}$$

where $M_L$ and $M_{L'}$ denote the matrices corresponding to $L$ and $L'$, respectively. Applying this identity directly yields the following crucial observation.

**Theorem 3.** *Let $g, g' \in \mathbb{F}_{p^n}[x]$ be two DO polynomials within the same extended-affine equivalence class and let $Y' \in \mathcal{D}_{g'} \cap \mathrm{GL}(n, \mathbb{F}_p)$. Then, there exist elements $A \in \mathrm{GL}(n, \mathbb{F}_p)$ and $Y \in \mathcal{D}_g \cap \mathrm{GL}(n, \mathbb{F}_p)$ such that $\mathcal{D}_{g'} Y'^{-1} = A^{-1} \cdot (\mathcal{D}_g Y^{-1}) \cdot A$.*

*More precisely, with the extended-affine equivalence transformation in Lemma 6, if $Y' = M_{g',\beta}$, we have $Y = M_{g,L(\beta)}$ and $A = M_{L'}^{-1}$.*

**Remark 3.** Note that Theorem 3 is a well-known observation for spread sets of presemifields: If $\Sigma, \Sigma'$ are the spread sets of two isotopic presemifields of order $p^n$, then $\Sigma' = X^{-1} \cdot \Sigma \cdot Y$ for some $X, Y \in \mathrm{GL}(n, \mathbb{F}_p)$ ([11, Sec. 2.1]). We gave a formulation and proof of Theorem 3 for general DO polynomials $g, g'$.

Conjugating a matrix $H \in \mathrm{Mat}_{\mathbb{F}_p}(n, n)$ with some element $A \in \mathrm{GL}(n, \mathbb{F}_p)$ as $A^{-1} \cdot H \cdot A$ corresponds to choosing a different basis for representing the same linear mapping. This is known as the *similarity transformation* of a matrix $H$ and yields an equivalence relation on $\mathrm{Mat}_{\mathbb{F}_p}(n,n)$.

**Definition 1** *(See p. 419 in [12]).* Two matrices $H, H' \in \mathrm{Mat}_{\mathbb{F}_p}(n,n)$ are called *similar*, written $H \sim H'$ if there exists an element $A \in \mathrm{GL}(n, \mathbb{F}_p)$ such that $H' = A^{-1} \cdot H \cdot A$. We denote by $[H]_\sim$ the equivalence class of $H$ with respect to $\sim$.

### *4.1. An invariant for extended-affine equivalence of DO polynomials*

Using matrix similarity, Theorem 3 yields an invariant for the extended-affine equivalence relation of DO polynomials. Indeed, Theorem 3 implies that, for two DO polynomials $g, g' \in \mathbb{F}_{p^n}[x]$ within the same extended-affine equivalence class, the set of quotients in $\mathcal{D}_g$ is the same as the set of quotients in $\mathcal{D}_{g'}$ up to applying the *same* similarity transformation to its elements. In other words, if $g$ and $g'$ are extended-affine equivalent DO polynomials, we have the equality[4]

$$\{([X]_\sim, \mathrm{mult}_{[X]_\sim}) \mid X \in \mathrm{Quot}(\mathcal{D}_g)\} = \{([X]_\sim, \mathrm{mult}_{[X]_\sim}) \mid X \in \mathrm{Quot}(\mathcal{D}_{g'})\}.$$

A canonical (efficiently computable) representative of a similarity equivalence class can be given by the *rational canonical form* (sometimes called *Frobenius normal form*), see. e.g., [12, page 475] for a definition. Thus, by denoting the rational canonical form of a matrix $X$ by $\mathrm{rcf}(X)$, we have the following corollary.

**Corollary 2.** *For two extended-affine equivalent DO polynomials $g, g' \in \mathbb{F}_{p^n}[x]$, we have the equality* $\mathrm{rcf}(\mathrm{Quot}(\mathcal{D}_g)) = \mathrm{rcf}(\mathrm{Quot}(\mathcal{D}_{g'}))$, *where* $\mathrm{rcf}(\mathrm{Quot}(\mathcal{D}_g))$ *denotes the set* $\{(\mathrm{rcf}(X), \mathrm{mult}_{\mathrm{rcf}(X)}) \mid X \in \mathrm{Quot}(\mathcal{D}_g)\}$.

Note that the invariant described in Corollary 2 is of no use in case that $\mathcal{D}_g$ does not contain any invertible element.

## 5. The set of quotients in $\mathcal{D}_{x^2}$

One can characterize the equivalence class of the planar monomial $g(x) = x^2$ by the algebraic structure of $\mathrm{Quot}(\mathcal{D}_g)$. Since it is well known in the context of finite (pre)semifields, we omit the proof.

**Lemma 7.** *Let $p$ be an odd prime, $n$ a positive integer, and $g \in \mathbb{F}_{p^n}[x]$ be a DO polynomial. The following assertions are equivalent.*

(1) *$g$ is equivalent to $x^2 \in \mathbb{F}_{p^n}[x]$.*

[4] Several elements in $\mathrm{Quot}(\mathcal{D}_g)$ can be in the same similarity equivalence class $[X]_\sim$, so we consider the multiplicities, denoted $\mathrm{mult}_{[X]_\sim}$ as well.

(2) *All elements in* $\mathcal{D}_g \setminus \{0\}$ *are invertible and we have* $\mathcal{D}_g Y^{-1} = \mathcal{D}_g Z^{-1}$ *for all* $Y, Z \in \mathcal{D}_g \setminus \{0\}$. *In particular* $|\mathrm{Quot}(\mathcal{D}_g)| = p^n$.
(3) *The matrix* $M_{g,1}$ *is invertible and the set* $\mathcal{D}_g M_{g,1}^{-1}$ *together with the addition and multiplication of matrices is a field isomorphic to* $\mathbb{F}_{p^n}$.

**Corollary 3.** *Let* $p$ *be an odd prime,* $g(x) \in \mathbb{F}_{p^n}[x]$ *be a DO polynomial and let* $Y := M_{g,1}$. *Let* $\alpha_1, \alpha_2, \ldots, \alpha_n$ *be an* $\mathbb{F}_p$*-basis of* $\mathbb{F}_{p^n}$. *Then,* $g$ *is equivalent to* $x^2$ *if and only if* $Y$ *is invertible and the matrix algebra*

$$\mathbb{F}_p[M_{g,\alpha_1}Y^{-1}, M_{g,\alpha_2}Y^{-1}, \ldots, M_{g,\alpha_n}Y^{-1}]$$

*is a field of order* $p^n$.

One way to decide whether $g$ is equivalent to $x^2$ is then to compute the elements $M_{g,\alpha_1}Y^{-1}, M_{g,\alpha_2}Y^{-1}, \ldots, M_{g,\alpha_n}Y^{-1}$, which can be done by evaluating the polynomial $g$ on $\mathcal{O}(n^2)$ many elements in $\mathbb{F}_{p^n}$, check the invertibility of those matrices, and apply Algorithm 2 to decide whether $\mathbb{F}_p[M_{g,\alpha_1}Y^{-1}, \ldots, M_{g,\alpha_n}Y^{-1}]$ is a field of extension degree $n$ over $\mathbb{F}_p$. Hence, with Algorithm 2, we can decide whether $g$ is equivalent to $x^2$ using $\mathcal{O}(n^7 \log(p))$ elementary field operations.

**Remark 4.** While Algorithm 2 solves a more general problem, deciding equivalence of $g$ to $x^2$ can also be done differently by means of the nuclei of the corresponding presemifields. Indeed, the commutative presemifield $\mathcal{R}_g$ is isotopic to a finite field if and only if $(x \star y) \star z = x \star (y \star z)$ holds for all $x, y, z \in \mathbb{F}_{p^n}$. If $\phi_\alpha$ denotes the evaluation map of $\Delta_{g,\alpha}$ in $\mathbb{F}_{p^n}$, this is equivalent to $\forall x, y, z \in \mathbb{F}_{p^n} : \phi_x(\phi_y(z)) = \phi_z(\phi_x(y))$, by using the commutativity of $\star$. It is sufficient to verify this identity on an $\mathbb{F}_p$-basis of $\mathbb{F}_{p^n}$ (i.e., $x, y, z$ are only basis elements), thus yielding an algorithm of polynomial complexity $\mathcal{O}(n^5)$ field operations (evaluating $\phi_x(y)$ is a matrix-vector operation, so in $\mathcal{O}(n^2)$ and there are $n^3$ combinations $x, y, z$ to check).

**Remark 5.** Other efficiently verifiable and non-trivial *necessary conditions* for the equivalence of a polynomial to $x^2$ have been known before. For instance, Budaghyan and Helleseth [4] showed that any DO polynomial $g$ equivalent to $x^2 \in \mathbb{F}_{p^n}[x]$ must contain a monomial of the form $x^{2p^i}$ with $0 \leq i < n$.

## 6. On the set of quotients in $\mathcal{D}_{x^{p^k+1}}$

We conclude by studying the set $\mathrm{Quot}(\mathcal{D}_g)$ for planar DO monomials in general, i.e., $g(x) = x^{p^\ell(p^k+1)} \in \mathbb{F}_{p^n}[x]$ with $p$ being an odd prime and $n/\gcd(k,n)$ being odd [9]. In the following, we assume, w.l.o.g., $\ell = 0$. We have already discussed the case of $k = 0 \mod n$ in the previous section. If $k \neq 0 \mod n$, those monomials correspond to Albert's commutative twisted fields. We show that for any DO polynomial $h \in \mathbb{F}_{p^n}[x]$ equivalent to a planar monomial, the set $\mathrm{Quot}(\mathcal{D}_h)$ always contains the finite field of order $p^n$. More precisely, we show the following.

**Theorem 2.** *Let $p$ be an odd prime and $n$ a positive integer. Let $g(x) \in \mathbb{F}_{p^n}[x]$ be a planar DO monomial. For any $\alpha, \beta \in \mathbb{F}_{p^n}^*$, the element $X := M_{g,\beta}M_{g,\alpha}^{-1} \in \mathrm{Quot}(\mathcal{D}_g)$ generates a field isomorphic to $\mathbb{F}_p(\alpha^{-1}\beta)$ viz. $\mathbb{F}_p[X]$, and $\mathbb{F}_p[X] \subseteq \mathrm{Quot}(\mathcal{D}_g)$.*

Let us denote by $\phi_\alpha \colon \mathbb{F}_{p^n} \to \mathbb{F}_{p^n}, x \mapsto \alpha x^{p^k} + \alpha^{p^k} x$ the evaluation map of the linearized derivative $\Delta_{x^{p^k+1},\alpha} \in \mathbb{F}_{p^n}[x]$. It is well known that $\phi_\alpha$ is invertible if and only if $n/\gcd(k,n)$ is odd (see [9]). We have the following for the inverse, which is a special case of Thm. 2.1 of [34]. It can also be proven by straightforward calculation of $\phi_\alpha^{-1}(\phi_\alpha(x))$.

**Lemma 8** *(Special case of Thm. 2.1 of [34]). Let $k$ be such that $n/\gcd(k,n)$ is odd. Let $d := n/\gcd(k,n)$. For $\alpha \in \mathbb{F}_{p^n}^*$, the inverse of $\phi_\alpha \colon x \mapsto \alpha x^{p^k} + \alpha^{p^k} x$ is given by*

$$\phi_\alpha^{-1} \colon x \mapsto \frac{\alpha}{2} \cdot \sum_{i=0}^{d-1} (-1)^i \alpha^{-(p^k+1)p^{ki}} x^{p^{ki}}.$$

The following lemma is immediate.

**Lemma 9.** *Let $k$ be such that $n/\gcd(k,n)$ is odd and let $\phi_\alpha \colon x \mapsto \alpha x^{p^k} + \alpha^{p^k} x$. For any $\alpha, \beta \in \mathbb{F}_{p^n}^*$, we have $\phi_\beta(\phi_\alpha^{-1}(x)) = (\beta^{p^k} - \alpha^{p^k-1}\beta) \cdot \phi_\alpha^{-1}(x) + \alpha^{-1}\beta x$.*

**Proof.** This follows directly from $(\phi_\beta - \alpha^{-1}\beta \cdot \phi_\alpha)(x) = (\beta^{p^k} - \alpha^{p^k-1}\beta)x$. □

The monomial $g(x) = x^{p^k+1}$ admits a non-trivial self equivalence via $g(x) = \gamma^{-(p^k+1)} \cdot g(\gamma x)$, where $\gamma$ is an arbitrary non-zero element of $\mathbb{F}_{p^n}$. Then, the following lemma directly follows from Lemma 6 and Eq. (2).

**Lemma 10.** *Let $k$ be such that $n/\gcd(k,n)$ is odd and let $\phi_\alpha \colon x \mapsto \alpha x^{p^k} + \alpha^{p^k} x$. For any $\alpha, \beta, \gamma \in \mathbb{F}_{p^n}$, $\alpha, \gamma \neq 0$, we have $\phi_\beta(\phi_\alpha^{-1}(x)) = \gamma^{-(p^k+1)} \cdot \phi_{\gamma\beta}(\phi_{\gamma\alpha}^{-1}(\gamma^{p^k+1} x))$.*

To show Theorem 2, we will first deduce that each element in $\mathrm{Quot}(\mathcal{D}_g)$ generates (a subfield of) $\mathbb{F}_{p^n}$. To do so, we show that each element in $\mathrm{Quot}(\mathcal{D}_g)$ is similar to a matrix corresponding to multiplication with an element of $\mathbb{F}_{p^n}$.

**Lemma 11.** *Let $k$ be such that $n/\gcd(k,n)$ is odd. Let $\alpha, \beta \in \mathbb{F}_{p^n}, \alpha \neq 0$. If $\alpha^{-1}\beta \in \mathbb{F}_{p^{\gcd(k,n)}}$, the mapping $\phi_\beta \circ \phi_\alpha^{-1}$ is equal to $x \mapsto \alpha^{-1}\beta x$. If $\alpha^{-1}\beta$ lies not in $\mathbb{F}_{p^{\gcd(k,n)}}$, the mapping $\psi_{\alpha,\beta} \circ \phi_\beta \circ \phi_\alpha^{-1} \circ \psi_{\alpha,\beta}^{-1}$ is equal to $x \mapsto (\alpha^{-1}\beta)^{p^k} x$, where*

$$\psi_{\alpha,\beta} \colon x \mapsto \alpha^{p^k} \cdot \phi_\alpha\left(\frac{1}{\beta^{p^k} - \alpha^{p^k-1}\beta} \cdot x\right).$$

**Proof.** We first observe that $\beta^{p^k} - \alpha^{p^k-1}\beta$ is equal to zero if and only if $\beta = 0$ or $(\alpha^{-1}\beta)^{p^k-1} = 1$, i.e., if and only if $\alpha^{-1}\beta$ is contained in the subfield $\mathbb{F}_{p^{\gcd(k,n)}} \subseteq \mathbb{F}_{p^n}$. Hence, by Lemma 9, the statement is trivial for the case of $\alpha^{-1}\beta \in \mathbb{F}_{p^{\gcd(k,n)}} \subseteq \mathbb{F}_{p^n}$.

In the other case, the mapping $\psi_{\alpha,\beta}$ is well defined and we can decompose $\psi_{\alpha,\beta}$ as $C \circ B \circ A$, where $A$ is a multiplication by $(\beta^{p^k} - \alpha^{p^k-1}\beta)^{-1}$, $B = \phi_\alpha$, and $C$ is a multiplication by $\alpha^{p^k}$. For all $x \in \mathbb{F}_{p^n}$, we then have:

$$L_1(x) := A(\phi_\beta(\phi_\alpha^{-1}(A^{-1}(x)))) = \phi_\alpha^{-1}\left((\beta^{p^k} - \alpha^{p^k-1}\beta)x\right) + \alpha^{-1}\beta x.$$

$$\begin{aligned} L_2(x) &:= B(L_1(B^{-1}(x))) = (\beta^{p^k} - \alpha^{p^k-1}\beta) \cdot \phi_\alpha^{-1}(x) + \phi_\alpha(\alpha^{-1}\beta \cdot \phi_\alpha^{-1}(x)) \\ &= \beta^{p^k} \cdot \left(\phi_\alpha^{-1}(x) + \alpha^{-p^k+1}(\phi_\alpha^{-1}(x))^{p^k}\right). \end{aligned}$$

$$\begin{aligned} L_3(x) &:= C(L_2(C^{-1}(x))) = \beta^{p^k} \cdot \left(\alpha^{p^k}\phi_\alpha^{-1}(\alpha^{-p^k}x) + \alpha(\phi_\alpha^{-1}(\alpha^{-p^k}x))^{p^k}\right) \\ &= \beta^{p^k} \cdot \phi_\alpha(\phi_\alpha^{-1}(\alpha^{-p^k}x)) = (\alpha^{-1}\beta)^{p^k}x. \end{aligned}$$

The proof is complete since $L_3 = \psi_{\alpha,\beta} \circ \phi_\beta \circ \phi_\alpha^{-1} \circ \psi_{\alpha,\beta}^{-1}$. □

**Remark 6.** In [2, Lem. 4], for the case of $k \mid n$ and $n/k$ odd, Albert showed that the characteristic polynomial of the mapping $\phi_\beta \circ \phi_1^{-1}$ is equal to the characteristic polynomial of the field multiplication by $\beta$, hence $\phi_\beta \circ \phi_1^{-1}$ generates a finite field (and then, so does $\phi_\beta \circ \phi_\alpha^{-1}$ by Lemma 10, choosing $\gamma = \alpha^{-1}$). More generally, let $\mathbb{F}$ be a finite extension of degree at least 3 over the field $\mathbb{F}_q$ with $q$ elements, $q > 2$. For $c \in \mathbb{F}^*$ not being a $(q-1)$-st power and $\beta \in \mathbb{F}$, let $\varphi_\beta^{(c)} : x \mapsto \beta^q x - c\beta x^q$ and $\widetilde{\varphi}_\beta^{(c)} : x \mapsto \beta x^q - c\beta^q x$, which are $\mathbb{F}_q$-linear mappings over $\mathbb{F}$. Albert showed that the characteristic polynomials of $\varphi_\beta^{(c)} \circ ({\varphi_1}^{(c)})^{-1}$ and $\widetilde{\varphi}_\beta^{(c)} \circ ({\widetilde{\varphi}_1}^{(c)})^{-1}$ are equal to the characteristic polynomial of the field multiplication by $\beta$.

**Remark 7.** The fact that an element $X := M_{g,\beta}M_{g,\alpha}^{-1}$ generates a finite field (which was already known before) implies that $X$ has an irreducible minimal polynomial. Hence, if $g(x) \in \mathbb{F}_{p^n}[x]$ is a DO polynomial equivalent to a planar monomial, each element of $\mathrm{Quot}(\mathcal{D}_g)$ has an irreducible minimal polynomial. We are not aware of planar DO polynomials inequivalent to monomials with this property, and we checked the known instances for $p = 3, n \leq 8$, listed in [16]. One could then check whether a given DO polynomial $g$ is inequivalent to a planar monomial by evaluating the minimal polynomials for several random elements in $\mathrm{Quot}(\mathcal{D}_g)$ and finding a reducible one. Of course, it does not yield a method to decide the equivalence, and we leave it an open question to provide a polynomial-time algorithm to do so.

Recall that an element $a$ in $\mathbb{F}_{p^n}$ is called a *square* in $\mathbb{F}_{p^n}$ if and only if the polynomial $x^2 - a$ splits in $\mathbb{F}_{p^n}$ and called a *non-square* otherwise. Moreover, it is well known that the set of non-zero squares $S$ in $\mathbb{F}_{p^n}$ is a subgroup of $\mathbb{F}_{p^n}^*$ of cardinality $\frac{p^n-1}{2}$. Thus, the factor group $\mathbb{F}_{p^n}^*/S$ consists of two elements and we have $\mathbb{F}_{p^n} = S \cup uS \cup \{0\}$, where $u$ is a non-square. Moreover it is well known that this result generalizes as follows: If, for an even value of $d$, the power mapping $x \mapsto x^d$ is 2-to-1 over $\mathbb{F}_{p^n}$, then the image of this mapping is $S \cup \{0\}$ and every square in $\mathbb{F}_{p^n}$ can be written as $\alpha^d$ for some $\alpha \in \mathbb{F}_{p^n}$, and every non-square can be written as $u\alpha^d$ for a properly chosen $\alpha \in \mathbb{F}_{p^n}$ and a fixed but arbitrary non-square $u \in \mathbb{F}_{p^n}$. The following well-known lemma gives a necessary condition when a non-square stays a non-square under a finite field extension. For the sake of completeness we include a proof.

**Lemma 12.** *Given a finite field extension $\mathbb{F}_{p^n}/\mathbb{F}_{p^m}$ of odd extension degree. If $a$ is a non-square in $\mathbb{F}_{p^m}$, then also in $\mathbb{F}_{p^n}$.*

**Proof.** If $a$ is a non-square in $\mathbb{F}_{p^m}$, the polynomial $x^2 - a$ is irreducible over $\mathbb{F}_{p^m}$. Hence, the splitting field is given by $\mathbb{F}_{p^m}(\sqrt{a})$, where $\sqrt{a}$ denotes one of the two roots of $x^2 - a$ in the algebraic closure, and we have $[\mathbb{F}_{p^m}(\sqrt{a}) : \mathbb{F}_{p^m}] = 2$. Assume now that $a$ is a square in $\mathbb{F}_{p^n}$. Then, $\mathbb{F}_{p^m}(\sqrt{a})$ is a subfield of $\mathbb{F}_{p^n}$ and hence the extension degree of $\mathbb{F}_{p^n}/\mathbb{F}_{p^m}$ is even, a contradiction. □

**Proof of Theorem 2.** Let $\alpha, \beta \in \mathbb{F}_{p^n}^*$ and let $X := M_{g,\beta}M_{g,\alpha}^{-1}$. By Lemma 11, the linear mapping $\phi_\beta \circ \phi_\alpha^{-1}$ is similar to $x \mapsto \alpha^{-1}\beta x$. Hence, the $\mathbb{F}_p$-algebra $\mathbb{F}_p[X]$ is isomorphic to $\mathbb{F}_p(\alpha^{-1}\beta)$ and thus a field. It is left to show that $\mathbb{F}_p[X] \subseteq \mathrm{Quot}(\mathcal{D}_g)$. The case of $\alpha^{-1}\beta \in \mathbb{F}_{p^{\gcd(k,n)}}$ is trivial and we therefore assume in the following that $\alpha^{-1}\beta \notin \mathbb{F}_{p^{\gcd(k,n)}}$. We will first handle the case of $\alpha = 1$ and show that $\left(M_{g,\beta}M_{g,1}^{-1}\right)^r \in \mathrm{Quot}(\mathcal{D}_g)$ for any integer $r \geq 2$. By Lemma 11, we have

$$\psi_{1,\beta} \circ \left(\phi_\beta \circ \phi_1^{-1}\right)^r \circ \psi_{1,\beta}^{-1}(x) = \left(\psi_{1,\beta} \circ \phi_\beta \circ \phi_1^{-1} \circ \psi_{1,\beta}^{-1}\right)^r (x) = \beta^{rp^k} x.$$

Further,

$$\beta^{rp^k} x = \begin{cases} \psi_{1,\beta^r} \circ \phi_{\beta^r} \circ \phi_1^{-1} \circ \psi_{1,\beta^r}^{-1}(x) & \text{if } \beta^r \notin \mathbb{F}_{p^{\gcd(k,n)}} \\ \beta^r x = \phi_{\beta^r} \circ \phi_1^{-1}(x) & \text{otherwise} \end{cases},$$

and thus

$$\left(\phi_\beta \circ \phi_1^{-1}\right)^r = \begin{cases} \psi_{1,\beta}^{-1} \circ \psi_{1,\beta^r} \circ \phi_{\beta^r} \circ \phi_1^{-1} \circ \psi_{1,\beta^r}^{-1} \circ \psi_{1,\beta} & \text{if } \beta^r \notin \mathbb{F}_{p^{\gcd(k,n)}} \\ \psi_{1,\beta}^{-1} \circ \phi_{\beta^r} \circ \phi_1^{-1} \circ \psi_{1,\beta} & \text{otherwise} \end{cases}. \tag{3}$$

We will now prove that the latter composition is equal to $\phi_\delta \circ \phi_\gamma^{-1}$ for properly chosen field elements $\delta, \gamma$.

*Case $\beta^r \in \mathbb{F}_{p^{\gcd(k,n)}}$.* In this case, $\left(\phi_\beta \circ \phi_1^{-1}\right)^r (x) = \psi_{1,\beta}^{-1} \circ \phi_{\beta^r} \circ \phi_1^{-1} \circ \psi_{1,\beta}(x) = \psi_{1,\beta}^{-1}(\beta^r \cdot \psi_{1,\beta}(x)) = \beta^r \cdot \psi_{1,\beta}^{-1}(\psi_{1,\beta}(x)) = \beta^r x = \phi_{\beta^r} \circ \phi_1^{-1}(x)$, since $\psi_{1,\beta}$ is $\mathbb{F}_{p^{\gcd(k,n)}}$-linear.

*Case* $\beta^r \notin \mathbb{F}_{p^{\gcd(k,n)}}$. We first observe that $\psi_{1,\beta}^{-1} \circ \psi_{1,\beta^r}(x) = \frac{\beta^{p^k}-\beta}{\beta^{rp^k}-\beta^r}x$. Let us define $\lambda := \frac{\beta^{p^k}-\beta}{\beta^{rp^k}-\beta^r} \in \mathbb{F}_{p^n}^*$. The image of the mapping $x \mapsto x^{p^k+1}$ over $\mathbb{F}_{p^n}$ is equal to the set of squares in $\mathbb{F}_{p^n}$. Indeed, every element in the image is a square as $p^k+1$ is even, and $x \mapsto x^{p^k+1}$ is 2-to-1 as a DO planar function [33]. Hence, if $\lambda$ is a square, we have $\lambda = \gamma^{p^k+1}$ for an element $\gamma \in \mathbb{F}_{p^n}^*$ and, otherwise, we have $\lambda = u\gamma^{p^k+1}$ with $u \in \mathbb{F}_{p^n}^*$ being an arbitrary non-square. Note that we can always choose $u \in \mathbb{F}_{p^{\gcd(k,n)}}^*$. Indeed, let $n = 2^m\ell$ and $k = 2^{m'}\ell'$ with $\ell, \ell'$ being odd, we necessarily have $m' \geq m$, as otherwise $n/\gcd(k,n)$ would be even. So, $\mathbb{F}_{p^{\gcd(k,n)}}$ contains $\mathbb{F}_{p^{2^m}}$ as a subfield and the extension degree $[\mathbb{F}_{p^n} : \mathbb{F}_{p^{\gcd(k,n)}}]$ is odd. The claim then follows by Lemma 12.

Let us therefore assume in the following that $\lambda = u\gamma^{p^k+1}$ with $\gamma \in \mathbb{F}_{p^n}^*$ and $u \in \mathbb{F}_{p^{\gcd(k,n)}}^*$. We have

$$\begin{aligned} \psi_{1,\beta}^{-1} \circ \psi_{1,\beta^r} \circ \phi_{\beta^r} \circ \phi_1^{-1} \circ \psi_{1,\beta^r}^{-1} \circ \psi_{1,\beta}(x) &= \lambda \cdot \left(\phi_{\beta^r} \circ \phi_1^{-1}\right)\left(\lambda^{-1}x\right) \\ &= \gamma^{p^k+1} \cdot \left(\phi_{\beta^r} \circ \phi_1^{-1}\right)\left(\gamma^{-(p^k+1)}x\right), \end{aligned} \tag{4}$$

where the last equality follows from the fact that $u \in \mathbb{F}_{p^{\gcd(k,n)}}^*$. By Lemma 10, we have $\gamma^{p^k+1} \cdot \left(\phi_{\beta^r} \circ \phi_1^{-1}\right)\left(\gamma^{-(p^k+1)}x\right) = \phi_{\gamma\beta^r} \circ \phi_\gamma^{-1}(x)$.

To handle the case of $\alpha \neq 1$, we apply Lemma 10 with $\gamma = \alpha^{-1}$ and obtain $\phi_\beta(\phi_\alpha^{-1}(x)) = \alpha^{p^k+1} \cdot \phi_{\alpha^{-1}\beta}(\phi_1^{-1}(\alpha^{-(p^k+1)}x))$, hence,

$$\begin{aligned} (\phi_\beta \circ \phi_\alpha^{-1})^r(x) &= \alpha^{p^k+1} \cdot (\phi_{\alpha^{-1}\beta} \circ \phi_1^{-1})^r(\alpha^{-(p^k+1)}x) \\ &= \alpha^{p^k+1} \cdot \left(\phi_{\delta'} \circ \phi_{\gamma'}^{-1}(\alpha^{-(p^k+1)}x)\right) = \phi_{\alpha\delta'} \circ \phi_{\alpha\gamma'}^{-1}(x) \end{aligned}$$

for appropriate elements $\gamma', \delta'$. We have now established that, for $\alpha^{-1}\beta$ being a generator of $\mathbb{F}_{p^n}^*$, the algebra $\mathbb{F}_p[X]$ is a field of order $p^n$ contained in $\mathrm{Quot}(\mathcal{D}_g)$.

To handle the general case where $\alpha^{-1}\beta$ is not a generator of $\mathbb{F}_{p^n}^*$, we will show that $X$ is equal to $(M_{g,\beta'}M_{g,\alpha'}^{-1})^r$ for some generator $\alpha'^{-1}\beta'$ of $\mathbb{F}_{p^n}^*$ and some non-negative integer $r$. Then, it would immediately follow that $\mathbb{F}_p[X] \subseteq \mathbb{F}_p[M_{g,\beta'}M_{g,\alpha'}^{-1}] \subseteq \mathrm{Quot}(\mathcal{D}_g)$. Indeed, let $\bar{\beta}$ be a generator of $\mathbb{F}_{p^n}^*$ such that $\bar{\beta}^r = \alpha^{-1}\beta$ and let

$$\frac{\bar{\beta}^{p^k} - \bar{\beta}}{\bar{\beta}^{rp^k} - \bar{\beta}^r} = u\gamma^{p^k+1}$$

with $\gamma \in \mathbb{F}_{p^n}^*$ and $u \in \mathbb{F}_{p^{\gcd(k,n)}}^*$. By extensively applying Lemma 10 and the result we established above, we obtain

$$\begin{aligned} (\phi_{\alpha\gamma^{-1}\bar{\beta}} \circ \phi_{\alpha\gamma^{-1}}^{-1})^r(x) &= \left((\alpha^{-1}\gamma)^{-(p^k+1)} \cdot \phi_{\bar{\beta}} \circ \phi_1^{-1}((\alpha^{-1}\gamma)^{p^k+1}x)\right)^r \\ &= (\alpha^{-1}\gamma)^{-(p^k+1)} \cdot \left(\phi_{\bar{\beta}} \circ \phi_1^{-1}\right)^r\left((\alpha^{-1}\gamma)^{p^k+1}x\right) \end{aligned}$$

$$= (\alpha^{-1}\gamma)^{-(p^k+1)} \cdot \phi_{\gamma\bar{\beta}^r} \circ \phi_{\gamma}^{-1}((\alpha^{-1}\gamma)^{p^k+1}x)$$
$$= (\alpha^{-1}\gamma)^{-(p^k+1)} \cdot \phi_{\alpha^{-1}\gamma\beta} \circ \phi_{\gamma}^{-1}((\alpha^{-1}\gamma)^{p^k+1}x)$$
$$= \phi_{\beta} \circ \phi_{\alpha}^{-1}(x). \quad \square$$

**Remark 8.** For a planar DO polynomial $g(x) = x^{p^k+1} \in \mathbb{F}_{p^n}[x]$, we can determine the cardinality of $\mathrm{Quot}(\mathcal{D}_g)$. To do so, we first show that, given $\beta, \bar{\beta} \in \mathbb{F}_{p^n}$ and $\alpha, \bar{\alpha} \in \mathbb{F}_{p^n}^*$, the equality $\phi_\beta \circ \phi_\alpha^{-1} = \phi_{\bar{\beta}} \circ \phi_{\bar{\alpha}}^{-1}$ holds if and only if one of the following two conditions is fulfilled:

(1) $\alpha^{-1}\beta \in \mathbb{F}_{p^{\gcd(k,n)}}$ and $\alpha^{-1}\beta = \bar{\alpha}^{-1}\bar{\beta}$,
(2) $\alpha^{-1}\beta \notin \mathbb{F}_{p^{\gcd(k,n)}}$ and $\alpha^{-1}\beta = \bar{\alpha}^{-1}\bar{\beta}$ and $\bar{\alpha} = c\alpha$ for $c \in \mathbb{F}_{p^{\gcd(k,n)}}^*$.

Indeed, we can easily deduce that Condition (1) or Condition (2) are sufficient by using Lemma 9 and the $\mathbb{F}_{p^{\gcd(k,n)}}$-linearity of the mappings $\phi_\gamma$, $\gamma \in \mathbb{F}_{p^n}$, respectively.

Let us now focus on showing the converse. From Lemma 9, the coefficient of $x$ of $\phi_\beta(\phi_\alpha^{-1}(x))$, interpreted as a polynomial in $\mathbb{F}_{p^n}[x]$, is equal to

$$\frac{(\alpha^{-1}\beta)^{p^k} + \alpha^{-1}\beta}{2}.$$

Hence, for $\phi_\beta \circ \phi_\alpha^{-1}$ and $\phi_{\bar{\beta}} \circ \phi_{\bar{\alpha}}^{-1}$ being equal, we necessarily have by comparing coefficients that $(\alpha^{-1}\beta - \bar{\alpha}^{-1}\bar{\beta})^{p^k} = -(\alpha^{-1}\beta - \bar{\alpha}^{-1}\bar{\beta})$. Since the kernel of $x \mapsto x^{p^k} + x$ over $\mathbb{F}_{p^n}$ is trivial, it follows that $\alpha^{-1}\beta = \bar{\alpha}^{-1}\bar{\beta}$. Further, we observe that the coefficient of $x^{p^k}$ of $\phi_\beta(\phi_\alpha^{-1}(x))$ is equal to

$$\frac{\alpha^{-p^{2k}+1}(\alpha^{-1}\beta - (\alpha^{-1}\beta)^{p^k})}{2}.$$

Hence, if $\phi_\beta \circ \phi_\alpha^{-1} = \phi_{\bar{\beta}} \circ \phi_{\bar{\alpha}}^{-1}$, by comparing the coefficients and using the fact that $\alpha^{-1}\beta = \bar{\alpha}^{-1}\bar{\beta}$, we obtain either $\alpha^{-1}\beta - (\alpha^{-1}\beta)^{p^k} = 0$ (i.e., $\alpha^{-1}\beta \in \mathbb{F}_{p^{\gcd(k,n)}}$) or $\alpha^{-p^{2k}+1} = \bar{\alpha}^{-p^{2k}+1}$. If $\alpha^{-1}\beta \notin \mathbb{F}_{p^{\gcd(k,n)}}$, by $\alpha^{p^n-1-p^{2k}+1} = \bar{\alpha}^{p^n-1-p^{2k}+1}$, it follows that $\bar{\alpha} = c\alpha$ for $c \in \mathbb{F}_{p^{\gcd(k,n)}}^*$. That is because $\gcd(p^n - p^{2k}, p^n - 1) = \gcd(p^{2k} - 1, p^n - 1) = p^{\gcd(2k,n)} - 1 = p^{\gcd(k,n)} - 1$, where the last equality holds because $n/\gcd(k,n)$ is odd. So, the mapping $x \mapsto x^{p^n - p^{2k}}$ over $\mathbb{F}_{p^n}$ is $(p^{\gcd(k,n)} - 1)$-to-1.

Having established the above conditions, we can now determine $|\mathrm{Quot}(\mathcal{D}_g)|$ as follows: For each of the $p^n - 1$ elements $\alpha \in \mathbb{F}_{p^n}^*$, there exist $p^n - p^{\gcd(k,n)}$ elements $\beta \in \mathbb{F}_{p^n}$ such that $\alpha^{-1}\beta \notin \mathbb{F}_{p^{\gcd(k,n)}}$ and for each such a pair $(\alpha, \beta)$, the polynomial $\phi_\beta \circ \phi_\alpha^{-1}$ can be written in $p^{\gcd(k,n)} - 1$ many ways by Condition (2). Indeed the condition states that $\alpha^{-1}\beta$ and $\bar{\alpha}^{-1}\bar{\beta}$ yield the same mapping if and only if $\bar{\alpha} = c\alpha$ and $\bar{\beta} = c\beta$ for $c \in \mathbb{F}_{p^{\gcd(k,n)}}^*$. This then yields

$$|\mathrm{Quot}(\mathcal{D}_g)| = \frac{(p^n - p^{\gcd(k,n)}) \cdot (p^n - 1)}{p^{\gcd(k,n)} - 1} + p^{\gcd(k,n)}.$$

## Acknowledgment

The authors thank Lukas Kölsch and an anonymous reviewer for pointing out the fact that deciding equivalence to $x^2$ can be done by verifying associativity of the corresponding commutative presemifield (see Remark 4). The authors further thank another anonymous reviewer for their helpful comments and for pointing us to the result established in [2, Lem.4].

The results of Section 6 have been presented at the 8th International Workshop on Boolean Functions and their Applications (BFA 2023).

The first author is funded by Deutsche Forschungsgemeinschaft (DFG) under Germany's Excellence Strategy - EXC 2092 CASA - 390781972.

## Data availability

No data was used for the research described in the article.